\documentclass[a4paper]{article}

\usepackage{graphicx}
\usepackage{latexsym}
\usepackage{amsmath}
\usepackage{amssymb}

\newenvironment{settheorem}{\begin{trivlist}\item[\hspace*{\labelsep}%
\textbf{Theorem of ZFC}]}{\end{trivlist}}

\newtheorem{postulate}{Postulate}

\newcommand{\brak}[1]{\ensuremath{\left[#1\right]}}

\newcommand{\paren}[1]{\ensuremath{\left(#1\right)}}

\newcommand{\pr}[1]{\Pr\!\brak{#1}}
\newcommand{\card}[1]{\mathop{\mathrm{Card}}\!\paren{#1}}
\newcommand{\powerset}[1]{\mathop{\mathcal{P}}\!\paren{#1}}

\newcommand{\reals}{\mathbb{R}}

\begin{document}

\title{Throwing Darts at a Ruler: Unpacking the Intuition Behind Freiling's Axiom of Symmetry}
\author{Anthony B.\ Morton}
\date{September 2013}
\maketitle

\section{Freiling's Darts: Denying the Continuum Hypothesis by Intuition}

It is well-known that the Continuum Hypothesis (CH) is independent of the other axioms of Zermelo-Fraenkel set theory, taken to include the Axiom of Choice (AC).
This formulation of set theory, conventionally denoted ZFC, is generally taken as the most well-justified foundation for mathematics in a logicist-realist framework.
On a realist view, the axioms of ZFC are taken to be intuitively justified by correspondence with reality.
This raises the question of whether an intuitive justification exists for CH as an additional axiom, or conversely whether it is more intuitive to deny CH.

Some years ago, Freiling \cite{f:aos} provided an intuitive argument for an axiom that, when appended to ZFC, is equivalent to the denial of CH.
Stated somewhat loosely, Freiling's \emph{Axiom of Symmetry} (AS) is the following:
\begin{postulate}[Axiom of Symmetry]
Given any function $f : \reals \mapsto \reals_{\aleph_0}$ that assigns to any real number $x$ a countable set of real numbers $f(x) = S_x \subset \reals$, there exist $x,y \in \reals$ such that $y \not\in f(x)$ and $x \not\in f(y)$.
\end{postulate}
The crucial aspect of this postulate is that it applies to \emph{any} function $f$ mapping individual real numbers to countable sets.
The reader should have little difficulty finding many examples of \emph{specific} functions $f$ for which AS is evidently true: it is the extension to \emph{all} such functions that takes AS beyond the obvious.

The intuitive argument for AS---with the interval $[0,1]$ in place of $\reals$, which introduces no particular difficulty as the two sets have the same cardinality---is as follows.

Consider the uniform probability distribution on $[0,1]$, and sample two independent random variables $X$ and $Y$ with this distribution.
(Freiling asks us to imagine two dart throwers who each aim randomly at the interval and act independently of one another.)
Let the respective sample points be $x$ and $y$.
Now consider the countable sets $S_x$ and $S_y$ that the function $f$ associates with $x$ and $y$.
Because $S_x$ is countable and $[0,1]$ is not, it is \emph{intuitively} practically certain that it does not contain the random sample $y \in [0,1]$.
Similarly, as $S_y$ is countable, it almost certainly does not contain $x$.
(Again intuitively, one may assert that if we have been unfortunate enough to pick one of the sporadic cases $x \in S_y$, or vice versa, we could just take new samples $x,y$ and repeat until with practical certainty we find values satisfying AS.)

The appeal to intuition is required in this argument because AS is equivalent in ZFC to the denial of CH (as will be shown shortly), and we have known since G\"{o}del that CH is consistent with ZFC.
This implies that AS has no proof as a theorem of ZFC.
Where Freiling's argument relies on intuition is its appeal to statements of probability: namely that $\pr{x \in S_y} = 0$ and $\pr{y \in S_x} = 0$ because $S_x$ and $S_y$ are countable, hence $x$ and $y$ satisfying AS can be found with `effective certainty'.
Ascribing probabilities to arbitrary subsets of $\reals$ is a fraught exercise in general, because rigorous probability statements about infinite sets require the sets to be \emph{measurable}, and in ZFC there are subsets of $\reals$ that fail to be measurable, whether according to the usual Lebesgue measure on $\reals$ or with respect to any other countably additive measure.
Freiling's argument manages to avoid this immediate problem, because countable additivity ensures that any countable set has Lebesgue measure zero.
But a leap of intuition \emph{is} necessary in order to proceed from the probability statement to a guarantee that $x \not\in S_y \wedge y \not\in S_x$ \emph{will} (often enough) be the case in practice.
This is particularly so when the sets $S_x$ and $S_y$ have `fractal' structure that forbids any intuitive attempt to reason directly about sampling events involving these sets.

Sets with such complex structure are involved in most of the interesting applications of AS; in particular the proof that in ZFC it implies the denial of CH.
This proof goes through as follows.
AC is equivalent to the proposition that every set can be well-ordered, so let $\prec$ denote a relation that well-orders $[0,1]$ and let $f$ be the function that assigns to each $x \in [0,1]$ the set of all numbers that precede $x$ in the well-ordering, together with $x$.
That is, $f(x) = \{x\} \cup \{t : t \in [0,1], t \prec x\}$.
It is a theorem of ZFC that any initial segment of a well-ordering on $[0,1]$, and therefore each $f(x)$ for arbitrary $x$, has cardinality strictly less than that of $[0,1]$.
Therefore, if CH is true and the only infinite cardinal less than $\card{[0,1]} = 2^{\aleph_0}$ is $\aleph_0$, then $f(x)$ for any $x$ has maximum cardinality $\aleph_0$ and is countable.
AS therefore applies and asserts that there exist $x,y$ with $y \succ x$ and $x \succ y$, a contradiction.

The argument for AS has attracted its fair share of objections, not least because in a sense it is \emph{too} persuasive.
If arbitrary countable sets satisfy AS because they have zero probability relative to a uniform distribution on a real interval, could the same intuition suggest the same for \emph{any} sets having cardinality less than $2^{\aleph_0}$, even if these may fail to be countable or indeed measurable?
Freiling's argument can easily be taken as suggesting an intuitive concept of `meagre probability' that relies not on strict measure theory, but simply on the cardinality of an `event' set relative to the `ambient' set.
It would therefore appear intuitive to assert the following stronger version of AS (call it SAS):
\begin{postulate}[Strong Axiom of Symmetry]
Given any function $f : \reals \mapsto \reals_{<2^{\aleph_0}}$ that assigns to any real number $x$ a set of real numbers $f(x) = S_x \subset \reals$ of cardinality less than that of $\reals$, there exist $x,y \in \reals$ such that $y \not\in f(x)$ and $x \not\in f(y)$.
\end{postulate}
Yet if one \emph{does} admit SAS, one can no longer assert AC!
After all, the same argument above that refutes CH will, given SAS, disallow \emph{any} well-ordering of $[0,1]$---thereby implying the negation of AC.
It appears that Freiling himself in \cite{f:aos} was persuaded by the ``case against AC'' implied by SAS (while ultimately settling for a version of SAS asserted only for sets $f(x)$ of Lebesgue measure zero).
Most theorists baulk at dropping AC from set theory given it provides a convenient framework for much of modern analysis (for example, by allowing the assertion that every vector space has a basis)---not to mention that strong intuitive arguments can be mounted for AC itself \cite{m:ba1,f:oaszfc}.

Commentators on Freiling such as Devlin \cite{d:hmrnat} and Brown \cite{b:pom} concede that the dart-throwing argument for AS is fundamentally intuitive and not rigorous.
Nonetheless, they are strongly persuaded by the argument---and by extension its refutation of CH---though they regard intuition alone as too weak to support SAS.
It is the intuitive content of Freiling's argument for AS, specifically the appeal to probabilities, that this note seeks to explore.

\section{The Inclusion Property: Unbelievable Monster, or Mere Counterexample?}

Abram \cite{a:npoch} treats AS from the converse standpoint: though the two statements are equivalent, one may choose to emphasise whether AS provides a \emph{denial} of CH, or conversely, whether CH provides a \emph{counterexample} to AS.
In Abram's view there is little to guide one's intuition apart from prior convictions about the validity of CH.
Given the relatively straightforward counterexample furnished by CH, it may be that one's intuition about AS itself is faulty:
\begin{quote}
That under the Continuum Hypothesis we may conceive of a specific function that negates the Axiom of Symmetry is a strong indication that the intuition of the axiom may be flawed (as intuition so commonly is).
Indeed, the negation of Freiling's Axiom of Symmetry is no more counterintuitive than that of a space-filling curve or a nowhere differentiable, everywhere continuous function, but these objects are known to exist.
\end{quote}

Indeed, the following is a straightforward corollary of the well-ordering principle and standard results about well-orderings.
\begin{settheorem}
Given any set $\Omega$, there exists a function $f : \Omega \mapsto \powerset{\Omega}$, assigning to each element of $\Omega$ a subset of $\Omega$ with the following properties:
\begin{enumerate}
\item\label{prop:inclusion}
For any $x, y \in \Omega$ distinct, either $f(x) \subset f(y)$ or $f(y) \subset f(x)$.
Also, $x \in f(x)$ for any $x$.
\item\label{prop:lowcard}
If $\Omega \subseteq \reals$, then $f(x)$ for any $x$ has cardinality less than $2^{\aleph_0}$.
\end{enumerate}
\end{settheorem}
The explicit construction of $f(x)$, as the initial segment to $x$ of a well-ordering $\prec$ on $\Omega$ (which exists by virtue of AC), was given above.
Property \ref{prop:inclusion} is called here the \emph{inclusion property}: it follows directly from the equivalent property for two initial segments of a well-ordering.
Property \ref{prop:lowcard} is a consequence of the fact that all initial segments have a unique order type: accordingly no initial segment of a well-ordering of $\reals$ can attain the order type of $\reals$ itself and so must in turn be of lesser cardinality than $\card{\reals} = 2^{\aleph_0}$.
This remains true if one restricts the well-ordering and its initial segments to any subset of $\reals$.

Abram's point above is well made: in postulating that there is a function $f$ with the inclusion property, there is no greater affront to intuition than implied by the notion of a well-ordering of $\reals$ (a firm consequence of AC).
But nor does it particularly offend one's intuition to affirm, accepting CH, that all of the subsets assigned by $f$ are in addition countable.
Then, however, it contradicts the intuitive justification for AS.

\section{The Problem is Independence}

But consider the statement of AS again: it asserts the existence of $x$ and $y$ such that $y \not\in f(x)$ \emph{and} $x \not\in f(y)$.
So it is asserting the truth of a \emph{joint} proposition for some choice of $x$ and $y$.
Yet, it seeks to assert this on probabilistic grounds: stated more explicitly, the argument is that the separate `events' $y \not\in f(x)$ and $x \not\in f(y)$ are each true with probability 1, hence the joint `event' $y \not\in f(x) \wedge x \not\in f(y)$ is also true with probability 1 owing to the `independent' choice of $x$ and $y$.

It may be, however, that AS is relying for its intuitive appeal on a mistaken notion of independence.
In probability calculus, the probability of a joint event is given as the product of the probability of one event alone, and another probability \emph{conditional} on that event:
\begin{equation}
\pr{a \wedge b} = \pr{a} \pr{b | a} = \pr{b} \pr{a | b}.
\label{eq:prodlaw}
\end{equation}
To assert that events $a$ and $b$ are independent is to assert that the conditional probabilities $\pr{a | b}$ and $\pr{b | a}$ do not differ from the unqualified probabilities, $\pr{a}$ and $\pr{b}$ respectively.
But if a conditional probability, say $\pr{b | a}$, does differ from $\pr{b}$, this does not necessarily mean that $a$ has any \emph{causal} influence on $b$: it means rather that there is some essential \emph{logical} connection between the two events.
If Mary places a red card and a black card in two separate sealed envelopes and John selects one at random, then a week later John opens his envelope and finds a red card, this does not \emph{cause} the card in Mary's remaining envelope to be black: yet without knowing the contents of any envelope one may justly state that $\pr{B} = 1/2$, $\pr{B | r} = 1$ and $\pr{B | b} = 0$, where $B$ is the event `the card in Mary's envelope is black', while $r, b$ are the events `the card in John's envelope is red (respectively, black)'.
The events $B$ and $b$ are not independent.

In the Freiling thought experiment, the sampling probabilities for $x$ and $y$ are indeed independent, by stipulation.
But does this independence automatically extend to the sets $f(x)$ and $f(y)$?
Here is the crucial importance of the inclusion property: it creates a logical connection between the sets $f(x)$ and $f(y)$ (one must include the other) even where there is none between the points $x$ and $y$ themselves.
A logical connection is likewise created between the events $y \not\in f(x)$ and $x \not\in f(y)$, so that for example
\begin{equation}
\pr{y \not\in f(x)} = 1, \qquad \text{yet} \qquad \pr{y \not\in f(x) \mathop{|} x \not\in f(y)} = 0.
\label{eq:princ}
\end{equation}
The formula on the right is in fact a purely logical statement expressed in probability calculus, a consequence of the inclusion property of $f$.
Given the truth of the proposition $x \not\in f(y)$, it must follow that $f(x) \not\subset f(y)$ (since $x \in f(x)$), and hence $f(y) \subset f(x)$ (since either $f(x) \subset f(y)$ or $f(y) \subset f(x)$ is true).
But now $y \in f(y) \subset f(x)$, and therefore $y \not\in f(x)$ is a false statement and must be assigned probability 0.

A similar line of reasoning establishes that $\pr{x \not\in f(y) \mathop{|} y \not\in f(x)} = 0$: the logical symmetry is preserved.
Notice again that the essential relationships are logical, not causal.
It is irrelevant which of $x$ and $y$ are sampled `first': the only question is whether $f(x)$ \emph{covers} $f(y)$ (in which case $y \not\in f(x)$ is not satisfied) or whether $f(x)$ \emph{is covered by} $f(y)$ (in which case $x \not\in f(y)$ is not satisfied).
One or other will be true in any case.
Meanwhile, (\ref{eq:prodlaw}) ensures that while $y \not\in f(x)$ and $x \not\in f(y)$ individually have probability 1, their conjunction has probability 0 when $f$ satisfies the inclusion property.

Once again, there is no necessary offence to intuition involved.
Just as dense sets can have an empty intersection, so events that are intuitively `almost certainly true' on their own can still be mutually exclusive if there is an underlying correlation.
Consider a scenario with $N$ envelopes distributed one by one to $N$ persons $P_1$, $P_2$, \ldots, $P_N$: each envelope contains a red card apart from one which has been selected at random to contain a black card.
If $N$ is large enough, the probability of the event $E_k$ : `person $P_k$ has a red card' can be asserted as close to 1 as desired absent knowing which envelope contains the black card; yet the conjunction of all the `almost true' events $E_k$ is demonstrably false.

So it may be with Freiling's Axiom: here there are only two events, yet one should not leap too readily from valid intuitions about events in isolation to reason about the joint event.

\bibliographystyle{unsrt}
\bibliography{bayes}

\end{document}